\documentclass[twocolumn]{autart}    

\usepackage{graphicx}          

\usepackage{mathrsfs}
\usepackage{tikz}
\usetikzlibrary{shapes, arrows, patterns}
\usepackage{xspace}

\makeatletter
\DeclareRobustCommand{\etc}{%
	\@ifnextchar{.}%
	{etc}%
	{etc.\@\xspace}%
}
\makeatother
\usepackage[tight,footnotesize]{subfigure}

\usepackage{booktabs}
\usepackage{tabularx}
\usepackage[table]{colortbl} 
\usepackage{multirow}
\usepackage{afterpage}
\usepackage{amsfonts}
\usepackage{amsmath,amsfonts,amssymb}

\newcommand{\bblue}[1]{\textcolor{black}{#1}}
 
\newcommand{\rredd}[1]{\textcolor{black}{#1}} 

\newcommand{\blue}[1]{\textcolor{black}{#1}}

\begin{document}

\begin{frontmatter}

\title{Non-indexability of the Stochastic Appointment Scheduling Problem\thanksref{footnoteinfo}} 

\thanks[footnoteinfo]{Corresponding author M.~Jafarnia-Jahromi.}

\author[ECE]{Mehdi Jafarnia-Jahromi}\ead{mjafarni@usc.edu},    
\author[ECE-CS]{Rahul Jain}\ead{rahul.jain@usc.edu},               

\address[ECE]{ECE Department, University of Southern California, Los Angeles, CA, USA}  
\address[ECE-CS]{ECE, ISE \& CS Department, University of Southern California, Los Angeles, CA, USA}             


\begin{abstract}                          
Consider a set of jobs with independent random service times to be scheduled on a single machine. The jobs can be surgeries in an operating room, patients' appointments in outpatient clinics, etc. The challenge is to determine the optimal sequence and appointment times of jobs to minimize some function of the server idle time and service start-time delay. We introduce a generalized objective function of delay and idle time, and consider $l_1$-type  and $l_2$-type  cost functions as special cases of interest. Determining an index-based policy for the optimal sequence in which to schedule jobs has been an open problem for many years. For example, it was conjectured that `least variance first' (LVF) policy is optimal for the $l_1$-type objective. This is known to be true for the case of two jobs with specific distributions. A key result in this paper is that the optimal sequencing problem is non-indexable, i.e., neither the variance, nor any other such index can be used to determine the optimal sequence in which to schedule jobs for $l_1$ and $l_2$-type objectives. We then show that given a sequence in which to schedule the jobs, sample average approximation yields a solution which is statistically consistent.
\end{abstract}

\end{frontmatter}

\section{Introduction}

Scheduling is an important aspect for efficient resource utilization, and there is a vast literature on the topic dating back several decades (see for example, \cite{Ba74,CoMaMi03,ErGoDe15,Pi12}). The problem considered in this paper is stochastic appointment scheduling which has various applications in scheduling of surgeries at operating rooms, appointments in outpatient clinics, cargo ships at seaports, etc.

\noindent\textit{Problem Statement:} Stochastic appointment scheduling problem (ASP) has a simple statement: Consider a finite set of $n$ jobs to be scheduled on a single machine. Job durations are random with known distributions. If a job completes before the appointment time of the subsequent job, the server will remain idle. Conversely, if it lasts beyond the allocated slot, the following job will be delayed. We need to determine optimal appointment times so that the expectation of a function of idle time and delay is minimized. Thus, the ASP addresses two important questions. First, given the sequence of jobs, what are the optimal appointment times? This is called the \textit{scheduling} problem. Second, in what sequence should the jobs be served? This is called the \textit{sequencing} problem. We note that while we call it the \textit{sequencing problem} for brevity, it really refers to the joint sequencing-scheduling problem of determining both the optimal sequence  as well as the appointment times.

Examples of the scheduling problem are cargo ships at seaports that must be given time to berth in the order in which they arrive. Examples of the sequencing problem are surgeries in a single operating room that must be scheduled the day before in an order that minimizes certain metrics such as wait times and idle time.

Intuitively speaking, scheduling jobs with more uncertain durations first may lead to delay propagation through the schedule. This intuition has motivated many researchers to prove optimality of `least variance first' (LVF) policy. However, efforts beyond the case of $n=2$ have not been fruitful for specific distributions such as uniform and exponential \cite{We90,Wa99}. Recent numerical work of \cite{MaMaZiKr18} has argued that LVF is not the best heuristic for the sequencing problem in the case that idle time and delay unit costs are not balanced. They introduce `newsvendor' index as another heuristic that outperforms variance. However, no proof of optimality is provided. This controversy raises an important open question that: \textit{Is there an index (a map from a random variable to the reals) that yields the optimal sequence?}
	
In this paper, both the \textit{scheduling} and \textit{sequencing} problems are addressed. In the \textit{sequencing} problem, we introduce an index (a map from a random variable to the reals) and prove that it is the only possible candidate to return the optimal sequence. This candidate index reduces to the `Newsvendor' index  and variance index for $l_1$ and $l_2$-type cost functions, respectively. However, by providing counterexamples for optimality of variance and newsvendor index, we show that the \textit{sequencing} problem is not indexable in general. Moreover, the candidate index illuminates that the heuristic sub-optimal indexing policy one might use depends on the objective function. In contrast to what has been used in the literature for a long time, variance is not the candidate index for $l_1$-type cost function. Instead, newsvendor should be used. This theoretical result confirms the numerical evidence of \cite{MaMaZiKr18} that newsvendor outperforms variance.
	
In the \textit{scheduling} problem, \blue{we prove that the $l_1$-type objective function introduced by \cite{We90} is convex and there exists a solution to the stochastic optimization problem.} \blue{Moreover, sample average approximation (SAA) can be used to approximate the optimization problem. We prove that SAA gives a solution that is statistically consistent.} These provide a feasible computational method to compute optimal appointment times given the sequence. To the best of our knowledge this result is new from two perspectives. First, the result is proved for a generalized objective function $c$ \blue{(as long as it is convex)} which includes previously considered objective function in the literature ($l_1$-type objective) as a special case of interest. Second, the assumptions for consistency of SAA is considerably relaxed compared to the standard literature of SAA. In fact we prove that SAA gives a consistent result as long as there exists a schedule with finite cost.

\blue{
The main contributions of this paper are:
\begin{itemize}
\item It is rigorously proved that there exists no index that yields the optimal sequence in the stochastic appointment scheduling problem (Theorem \ref{thm:no opt index}). However, a candidate index is introduced that yields a heuristic sub-optimal sequence. This candidate index reduces to newsvendor and variance for $l_1$-type and $l_2$-type objective functions, respectively.  
\item It is proved that the objective function of \cite{We90} for stochastic appointment scheduling problem is convex and there exists an optimal schedule (Proposition \ref{prop: convexity of g1} in Appendix A and Theorem \ref{thm:convexity}).
\item It is proved that for a fixed sequence of jobs, sample average approximation yields an approximate solution for the optimal schedule that is statistically consistent (Theorem \ref{thm:saa}).
\end{itemize}
}

The remainder of the paper is organized as follows. In Section \ref{sec:literature}, literature is reviewed. Section \ref{sec:problem} provides the problem formulation. In Section \ref{sec:seq}, the \textit{sequencing} problem is discussed. The \textit{scheduling} problem is addressed in Section \ref{sec:sch}. Section \ref{sec:numerical} provides numerical results followed by conclusions in Section \ref{sec:conclusion}.

\section{Literature Review}
\label{sec:literature}
\blue{Extensive application of appointment scheduling in transportation such as bus scheduling \cite{wu2019stochastic}, as well as healthcare applications \cite{ErDe13} has led to a vast body of literature on the topic. We note that the problem considered in this paper is a variant of a broader scheme of \textit{railway} scheduling \cite{herroelen2004construction,tian2014railway} where jobs are not allowed to start earlier than their schedule. Other variants of scheduling such as \textit{roadrunner} scheduling \cite{newbold1998project} are beyond the scope of this work.} Among all the related papers, we focus on the most relevant work and classify it into \textit{sequencing} and \textit{scheduling}. The interested reader is referred to \blue{\cite{HuKoBoHaBa12,GuDe08,CaVe03,AhJaKl17,cardoen2010operating,kuiper2017optimal,ChRo14,Ba74,CoMaMi03,ErGoDe15,Pi12,KuMa15,KuMa15practical,wu2019stochastic,ErDe13,tian2014railway}} and references therein for other aspects of the problem. 
\paragraph*{\bf Sequencing.}
The \textit{sequencing} problem we consider was first formulated by \cite{We90}. Intuitively speaking, jobs with less uncertainty should be placed first to avoid delay propagation throughout the schedule. Motivated by this intuition, a large body of literature suggested `least variance first' (LVF) rule as a sequencing policy (see \cite{We90,Wa99,MaRoZh14,DeViVo07,Qi16}). However, optimality of LVF rule is only proved for the case of two jobs $n=2$ for certain distributions \cite{We90,Wa99}. \cite{MaRoZh14} and \cite{GuDaJaJu16} tried to impose conditions under which LVF rule is optimal but the conditions are relatively restrictive and unlikely to hold in most scenarios of interest. A variant of the problem where jobs are allowed to start before scheduled appointments (no idle time is allowed) is studied by \cite{GuDaJaJu16,Ba14}. In particular, \cite{GuDaJaJu16} shows that LVF rule is optimal if there exists a dilation ordering for service durations. \cite{Gu07,BeDeErRoHu14} proved that if there exists a convex ordering for job durations, it is optimal to schedule smaller in convex order first for $n=2$. However, their efforts for $n > 2$ have not been fruitful. \cite{KoLeTeZh16} considered likelihood ratio as a measure of variability and obtained some insights into why smallest variability first may not be optimal. Based on the insights, they provided a counterexample for non-optimality of LVF rule in the case of $n=6$.

Besides the theoretical work, some papers have resorted to extensive simulation studies to investigate optimality of heuristics (see \cite{DeViVo07,KlRo96,Le03,MaDe06}). In particular, \cite{DeViVo07} numerically showed that LVF outperforms sequencing in increasing order of mean and coefficient of variation. However, \cite{MaSt12,MaMaZiKr18} argued that LVF is not the best sequencing policy especially when idle time and delay cost units are not balanced. Alternatively, \cite{MaMaZiKr18} proposed a `newsvendor' index and supported its better performance in simulations. No proof of optimality was provided.

\paragraph*{\bf Scheduling.}
The \textit{scheduling} problem is also intensively studied in the literature. The seminal work of \cite{Ba52} recommended to set appointment intervals equal to the average service time of each job. This approach was further persued by \cite{So66} and \cite{ChBa16}. However, letting job slots to be average service time can be near optimal only in the case that waiting cost is about 10\% to 50\% of the idle cost (see \cite{DeGu03}).

Starting with \cite{We90}, some papers modeled the problem using stochastic optimization to optimize on slot duration. He considered weighted sum of idle time and delay as the objective function and noticed that for the case of $n=2$, the problem is equivalent to the newsvendor problem. Based on that, he proposed a heuristic estimate of the job start times for $n>2$. This heuristic was extended by \cite{KeKlMa14} to general convex function of idle time and delay. \cite{Wa93} and \cite{Wa99} considered another objective function as the weighted sum of jobs' flow time (delay and service time) and server completion time and proved its convexity. Assuming that the job durations are exponentially distributed, he provided a set of nonlinear equations to derive the optimum slot durations. \cite{ViKuKeBh15} presented a lag order approximation by ignoring the effect of previous jobs past a certain point. 
\cite{KoLeTeZh13} adopted a robust method over all distributions with a given mean and covariance matrix of job durations. They computationally solved for 36 jobs and showed that their  solution is within 2\% of the approximate optimal solution given by \cite{DeGu03}. We refer the reader to \cite{BeQu11} and \cite{KaKo07} for approaches in  discrete time.

\section{Problem Formulation}
\label{sec:problem}

We start with providing a mathematical formulation of the problem. Let $\mathscr{X}^+$ be the space of nonnegative random variables and $\mathbf X = (X_1, \cdots, X_n)$ be a vector of independent random variables with components in $\mathscr{X}^+$ and known distributions denoting jobs durations to be served on a single server. \blue{Let $\sigma \cdot \mathbf X$ denote a permutation of the jobs in the order they are served, such that $\sigma(i)$ is the $i$th job receiving service.} Without loss of generality, we can assume that the first job starts at time zero, i.e., $s_1 = 0$. Let  $\mathbf s = (s_2, \cdots, s_n)$ be appointment times for job $2$ through $n$ in the order $\sigma$ and let $E_{\sigma(i)}$ be a random variable denoting the end time of job $i$ in this order (see Figure \ref{fig:cases}). Job $\sigma(i)$ may finish before or after scheduled start time of the subsequent job. In the case that $E_{\sigma(i)} \leq s_{i+1}$, job $\sigma(i+1)$ starts according to the schedule and server is idle between $E_{\sigma(i)}$ and $s_{i+1}$. In the case where $E_{\sigma(i)} > s_{i+1}$, job $\sigma(i+1)$ is delayed by $E_{\sigma(i)} - s_{i+1}$ and will start as soon as the previous job is finished. Hence,
\begin{align}
&E_{\sigma(1)} = X_{\sigma(1)}, \nonumber \\
&E_{\sigma(i)} = \max (E_{\sigma(i-1)}, s_i) + X_{\sigma(i)}, \quad i=2, \cdots, n.
\label{eq:Ei}
\end{align}

\begin{figure}
	\begin{center}
		\scalebox{.8}{\begin{tikzpicture}
\draw[very thick, ->] (0, 0) -- (10, 0);

\begin{scope}
\draw [dashed] (0, 0) node[anchor=north]{$s_1=0$} -- (0, 1);
\draw [thick, <-] (0, 0.5) -- (.5, 0.5) node[anchor=west]{$X_{\sigma(1)}$};
\draw [thick, ->] (1.6, .5) -- (2, 0.5);
\draw [dashed] (2, 0) node[anchor=north]{$E_{\sigma(1)}$} -- (2, 0.5);
\end{scope}

\begin{scope}[shift={(3, 0)}]
\draw [dashed] (0, 0) node[anchor=north]{$s_2$} -- (0, 1);
\draw [thick, <-] (0, 0.5) -- (1.5, 0.5) node[anchor=west]{$X_{\sigma(2)}$};
\draw [thick, ->] (2.6, .5) -- (4, 0.5);
\draw [dashed] (4, 0) node[anchor=north]{$E_{\sigma(2)}$} -- (4, 0.5);
\end{scope}

\begin{scope}[shift={(6, 0)}]
\draw [dashed] (0, 0) node[anchor=north]{$s_3$} -- (0, 1);
\draw [thick, <-] (1, 0.5) -- (1.5, 0.5) node[anchor=west]{$X_{\sigma(3)}$};
\draw [thick, ->] (2.6, .5) -- (3, 0.5);
\draw [dashed] (3, 0) node[anchor=north]{$E_{\sigma(3)}$} -- (3, 0.5);
\end{scope}

\end{tikzpicture}}
		\caption{Appointment scheduling. $s_i$ denotes appointment time of job $i$. \bblue{For the realization shown in this figure, server remains idle between $E_{\sigma(1)}$ and $s_2$ and third job is delayed for $E_{\sigma(2)} - s_3$ amount of time.}}
		\label{fig:cases}
	\end{center}
\end{figure}
Our goal is to determine appointment times such that a combination of both delay and idle time is optimized. Consider an objective function of form,
\begin{align}
\label{eq: C objective function}
C(\mathbf s, \sigma \cdot \mathbf X) = \sum_{i=2}^ng(E_{\sigma(i-1)} - s_i)
\end{align}
where $g:\mathbb{R} \to \mathbb{R}$ is a nonnegative, continuous and coercive function (i.e., $\lim_{| t | \to \infty}g(t)=\infty$). Furthermore, we can assume that $g(0) = 0$ since a perfect scenario where $E_{\sigma(i-1)} = s_i$ should not impose any cost. However, this assumption is not technically necessary. \blue{For the special case of $g = g_1$ in Example \ref{ex:c1}, this objective function reduces to that of \cite{We90}. However, \eqref{eq: C objective function} is not the most general objective function one can consider. For example, in some applications, it is useful to distinguish between jobs by considering different delay per unit costs for different jobs. Moreover, \eqref{eq: C objective function} does not account for overtime, a related quantity that is important in some applications.}

Given the schedule $\mathbf s$, $C(\mathbf s, \sigma \cdot \mathbf X)$ captures the associative cost of the realization of job durations $\mathbf X$ in the order $\sigma$. Thus, $c^\sigma(\mathbf s) = \mathbb{E}[C(\mathbf s, \sigma \cdot \mathbf X)]$ denotes the expected cost of schedule $\mathbf s$ when the jobs are served in the order $\sigma$. 

In the scheduling problem in Section \ref{sec:sch}, we assume that the sequence of jobs is given, and we are looking for a schedule that minimizes the expected cost, i.e.,
\begin{align}
\label{eq:obj}
\inf_{\mathbf s \in \mathcal{S}} c^\sigma(\mathbf s)
\end{align}
where $\mathcal{S} = \{(s_2, \cdots, s_n) \in \mathbb{R}^{n-1} \mid 0\leq s_2 \leq \cdots \leq s_n\}$ is a closed and convex subset of $\mathbb{R}^{n-1}$. The sequencing problem discussed in Section \ref{sec:seq}, addresses the question of finding the optimal order and appointment times of the jobs, \blue{i.e.,
\begin{align*}
\min_\sigma \inf_{\mathbf s \in \mathcal{S}} c^\sigma(\mathbf s).
\end{align*}
}

Before proceeding with the sequencing problem, let's see some possible choices for the function $g$.

\begin{figure}
	
	\centering
	\includegraphics[scale=.4]{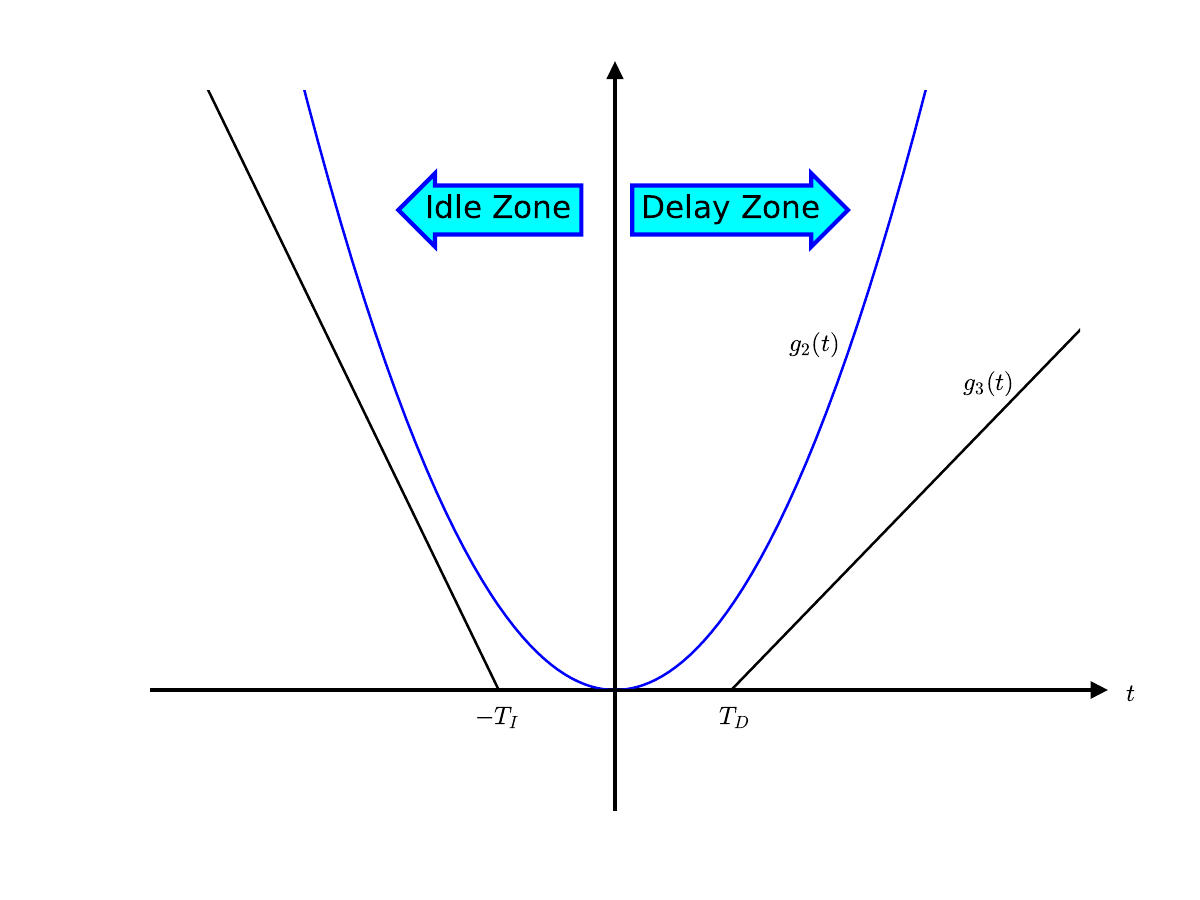}
	
	\caption{Examples of function $g$.}
	\label{fig:g}
	
\end{figure}

\begin{exmp}
	\normalfont
	\label{ex:c1}
	Let $g_1(t) = \beta(t)^+ + \alpha(-t)^+$ where $(\cdot)^+ = \max(\cdot, 0)$ and $\alpha , \beta > 0$. Thus, the objective function would be
	\begin{align}
	c^\sigma_1(\mathbf s) = \sum_{i=2}^n \mathbb{E}[\alpha (s_i - E_{\sigma(i-1)})^+ + \beta (E_{\sigma(i-1)} - s_i)^+].
	\end{align}
	$(s_i - E_{\sigma(i-1)})^+$ denotes idle time before job $i$ and $(E_{\sigma(i-1)} - s_i)^+$ indicates its possible delay. Cost function $c^\sigma_1$ is the same cost function used by \cite{We90}. If $\alpha \neq \beta$, it captures potential different costs associated with idle time and delay. We call this function $l_1$-type objective function.
\end{exmp}

\begin{exmp}
	\normalfont
	\label{ex:c2}
	Let $g_2(t) = t^2$. The objective function reduces to
	\begin{align}
	c^\sigma_2(\mathbf s) = \sum_{i=2}^n \mathbb{E}[(E_{\sigma(i-1)} - s_i)^2].
	\end{align}
	Cost function $c^\sigma_2$ penalizes both idle time and delay equally. However, due to the nonlinearity of $c^\sigma_2$, long idle time and delay are less tolerable. We call this function $l_2$-type objective function.
\end{exmp}

\begin{exmp}
	\normalfont
	\label{ex:g3}
	Let
	\begin{align}
	g_3(t) = \begin{cases}
	\beta (t - T_D), & \text{if } t \geq T_D  \\
	-\alpha (t + T_I), & \text{if } t \leq -T_I \\
	0, & \text{otherwise}
	\end{cases}
	\end{align}
	where $T_D, T_I \geq 0$ are delay and idle time tolerance, respectively (see Figure \ref{fig:g}). In this case, no cost is exposed for delay and idle time under a certain threshold. This situation arises in some applications such as operating room scheduling where some small amount of delay is tolerable.
\end{exmp}

\section{Sequencing Problem}\label{sec:seq}

\subsection{Non-indexability}

\bblue{In this section, we first consider the joint sequencing-scheduling problem (referred to as just the `sequencing problem' since the optimal sequence cannot be determined without also determining the optimal appointment times). } Intuitively, scheduling jobs with higher uncertainty in durations first may lead to delay propagation through the schedule. Considering objective function $c^\sigma_1$, this intuition has motivated many researchers to prove optimality of least variance first (LVF) policy. However, the efforts have not been fruitful beyond the case of two jobs ($n=2$) for some typical distributions such as exponential and uniform. Most related papers thus have resorted to numerical evaluation to analyze the performance of the LVF rule. In particular, \cite{DeViVo07} compared three ordering policies, namely, increasing mean, increasing variance, and increasing coefficient of variation. Using numerical experiment with real surgery duration data, they argued that ordering with increasing variance outperforms the other two heuristics. However, \cite{MaMaZiKr18} claimed that variance does not distinguish the potential difference between idle time and delay for $c^\sigma_1$. They introduced the \textit{newsvendor} index defined as
\begin{align}
\label{eq:optindex1}
I_1^*(X) =  \alpha \mathbb{E}[(s^* - X)^+] + \beta\mathbb{E}[(X - s^*)^+]
\end{align}
\bblue{where $F_X$ is cumulative distribution function of $X$}, $s^* := F_X^{-1}(\frac{\beta}{\alpha + \beta})$, and \bblue{numerically verified that sequencing in increasing order of $I_1^*$ outperforms LVF, and conjectured that it returns the optimal sequence. No proof of optimality was given. }
These conjectures will be evaluated in this section. In particular, we will prove that there exists no index (a map from a random variable to the reals) that yields the optimal sequence for objective functions $c^\sigma_1$ and $c^\sigma_2$.

\rredd{Moreover, we rigorously prove that the only candidate to provide the optimal sequence is newsvendor index for objective function $c^\sigma_1$ and variance for objective function $c^\sigma_2$. This provides a theoretical support for numerical evidence of \cite{MaMaZiKr18}. Moreover, it completely eliminates variance as a candidate heuristic for objective function $c^\sigma_1$.} 

Let's first start with a simple example of sequencing two jobs.
\begin{exmp}
	\label{ex:index motivation}
	\normalfont
	Consider the case of scheduling two jobs with durations $X_1, X_2$. The optimization problem to determine optimal appointment times given the sequence $(X_1, X_2)$ would be:
	\begin{align}
	\label{eq:tmp opt index}
	\inf_{s_2 \geq 0} \mathbb{E}[g(X_1 - s_2)]
	\end{align}
	The optimal cost given by the above equation is indeed an index that maps random variable $X_1$ to a real number. Moreover, sorting in increasing order of this index yields the optimal sequence for $n=2$.
\end{exmp}
Motivated by this example, we have a candidate index for general $n$:
\begin{align}
I_g^*(X) = \inf_{s \geq 0} \mathbb{E}[g(X - s)].
\end{align}
One can verify that this index reduces to variance ($I_2^*$) and newsvendor index ($I_1^*$) in the case that $g(t) = t^2$ and $g(t) = \beta (t)^+ + \alpha (-t)^+$, respectively. The natural question is whether this index provides the optimal sequence for $n > 2$. And if not, whether there is any other index that yields the optimal sequence. \bblue{In the ensuing, we will show that the answer to both of these questions is negative. In fact, we first prove in Proposition \ref{prop:optimal indexing} that $I_g^*$ is the only possible candidate to return the optimal sequence and then through counterexamples \ref{ex:l1} and \ref{ex:l2} show that it is not optimal.}

\bblue{To prepare the setup for Proposition \ref{prop:optimal indexing},} let $\bar{\mathbb{R}} = \mathbb{R} \cup \{+\infty\}$ be the extended real line. We say $I: \mathscr{X}^+ \to \bar{\mathbb{R}}$ is an index and denote the space of all indexes by $\mathscr{I}$. For example, mean, variance, newsvendor and $I_g^*$ are examples of elements in $\mathscr{I}$. First, we define an equivalence relation on $\mathscr{I}$.
\begin{defn}
	Let $I_1, I_2 \in \mathscr{I}$. We say $I_1$ is in relation with $I_2$ denoting by $I_1RI_2$ if for any $X_1, X_2 \in \mathscr{X}^+$, $I_1(X_1) \leq I_1(X_2)$ if and only if $I_2(X_1) \leq I_2(X_2)$.
\end{defn}
It is straightforward to check that $R$ is an equivalence relation on $\mathscr{I}$. Hence, $R$ splits $\mathscr{I}$ into disjoint equivalence classes. Next, we define a notation for sorting random variables in increasing order of an index.
\begin{defn}
	Let  $\mathbf X = (X_1, \cdots, X_n)$ be a random vector where $X_i \in \mathscr{X}^+$ for all $i$, and $I \in \mathscr{I}$ be an index. We say $\sigma \cdot \mathbf X = (X_{\sigma(1)}, \cdots, X_{\sigma(n)})$ is a valid permutation of $\mathbf X$ with respect to $I$ if $I(X_{\sigma(1)}) \leq \cdots \leq I(X_{\sigma(n)})$. We denote the set of valid permutations by $\mathscr{P}_I(\mathbf X)$.
	
	In the case that $I(X_1), \cdots, I(X_n)$ take distinct values, $\mathscr{P}_I(X)$ includes only one element.
\end{defn}
If $I_1$ is equivalent to $I_2$, then $\mathscr{P}_{I_1}(\mathbf X) = \mathscr{P}_{I_2}(\mathbf X)$ for any random vector $\mathbf X$ with components in $\mathscr{X}^+$.

\begin{defn}
	Index $I$ is optimal for cost function $c$ if for any $n \geq 2$ and any random vector $\mathbf X = (X_1, \cdots, X_n)$ with components in $\mathscr{X}^+$, $\inf_\mathbf s \mathbb{E}[C(\mathbf s,\sigma \cdot \mathbf X)] \leq \inf_\mathbf s \mathbb{E}[C(\mathbf s,\mathbf X)]$ for all $\sigma \in \mathscr{P}_I(\mathbf X)$.
\end{defn}
Thus, by the above remark if an index of a class is optimal, all equivalent indices are also optimal. Hence, optimality is a class property.

We already observed that $I_g^*$ is optimal for the case of $n=2$. The following Proposition provides a result for general $n$.
\begin{prop}
	\label{prop:optimal indexing}
	If there exists an optimal index for cost function $c$, it is equivalent to $I_g^*$.
\end{prop}

\begin{pf}
	Assume by contradiction that there exists index $J$ which is optimal but not equivalent to $I_g^*$. Hence, there exist random variables $X_1, X_2 \in \mathscr{X}^+$ such that $I_g^*(X_1) < I_g^*(X_2)$ but $J(X_1) \geq J(X_2)$. Note that $I_g^*(X_1) = \inf_{s \geq 0} \mathbb{E}[g(X_1 - s)]$ and $I_g^*(X_2) = \inf_{s \geq 0}\mathbb{E}[g(X_2 - s)]$. Hence, $I_g^*(X_1) < I_g^*(X_2)$ implies that $\inf_{s \geq 0} \mathbb{E}[g(X_1 - s)] < \inf_{s \geq 0}\mathbb{E}[g(X_2 - s)]$. However, optimality of $J$ implies that $\inf_{s \geq 0} \mathbb{E}[g(X_1 - s)] \geq \inf_{s \geq 0}\mathbb{E}[g(X_2 - s)]$ which is a contradiction. 
\qed
\end{pf}

Note that $I_g^*$ reduces to $I_1^*$ and $I_2^*$ for objective functions $c^\sigma_1$ and $c^\sigma_2$, respectively. We also notice that contrary to widely believed conjectures in the literature that $I_2^*$ (LVF rule) is an optimal index-based policy for cost function $c^\sigma_1$, Proposition \ref{prop:optimal indexing} states that variance can only be a candidate for $c^\sigma_2$. However, note that this proposition doesn't say anything about the existence of an optimal index. In the following, we provide counter examples which show that sequencing (and optimally scheduling) in increasing order of $I_1^*$ and $I_2^*$ is not optimal for $c^\sigma_1$ and $c^\sigma_2$, respectively.

\begin{exmp}
	\label{ex:l1}
	\normalfont
	Let $X_1, X_2, X_3$ be independent random variables in $L^1$ and assume that $X_1 \sim U(0,1)$ and $X_2$ follows the following distribution (see Figure \ref{fig:ex}):
	\[
	F_{X_2}(x)=\left\{
	\begin{array}{ll}
	0,& \text{if } x \leq 0\\
	2x^2,& \text{if } 0<x<0.5 \\
	2(x - 0.5)^2 + 0.5,& \text{if } 0.5\leq x<1 \\
	1,& \text{otherwise}
	\end{array}
	\right.
	\]
	
	\begin{figure}
		\centering     
		\subfigure[]{\label{fig:a}\includegraphics[width=80mm]{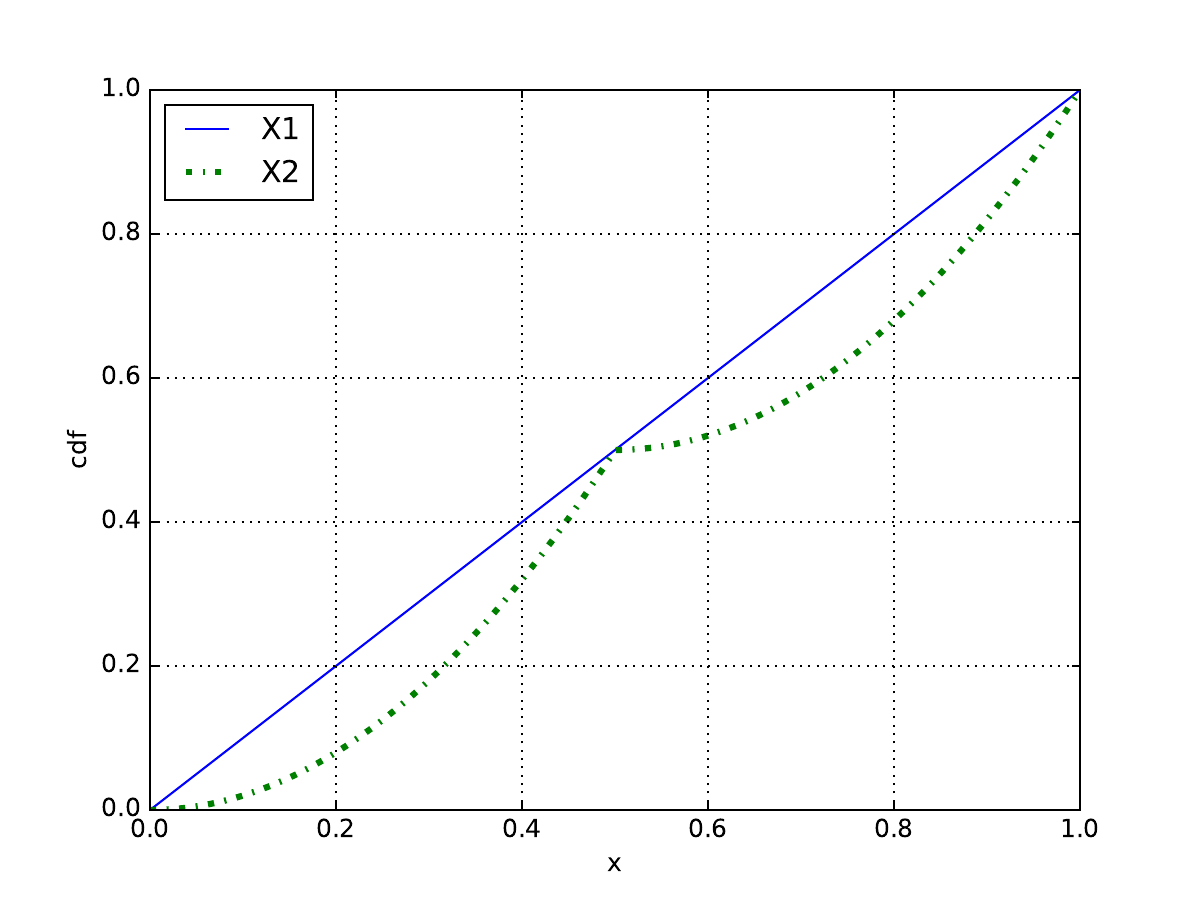}}
		\subfigure[]{\label{fig:b}\includegraphics[width=80mm]{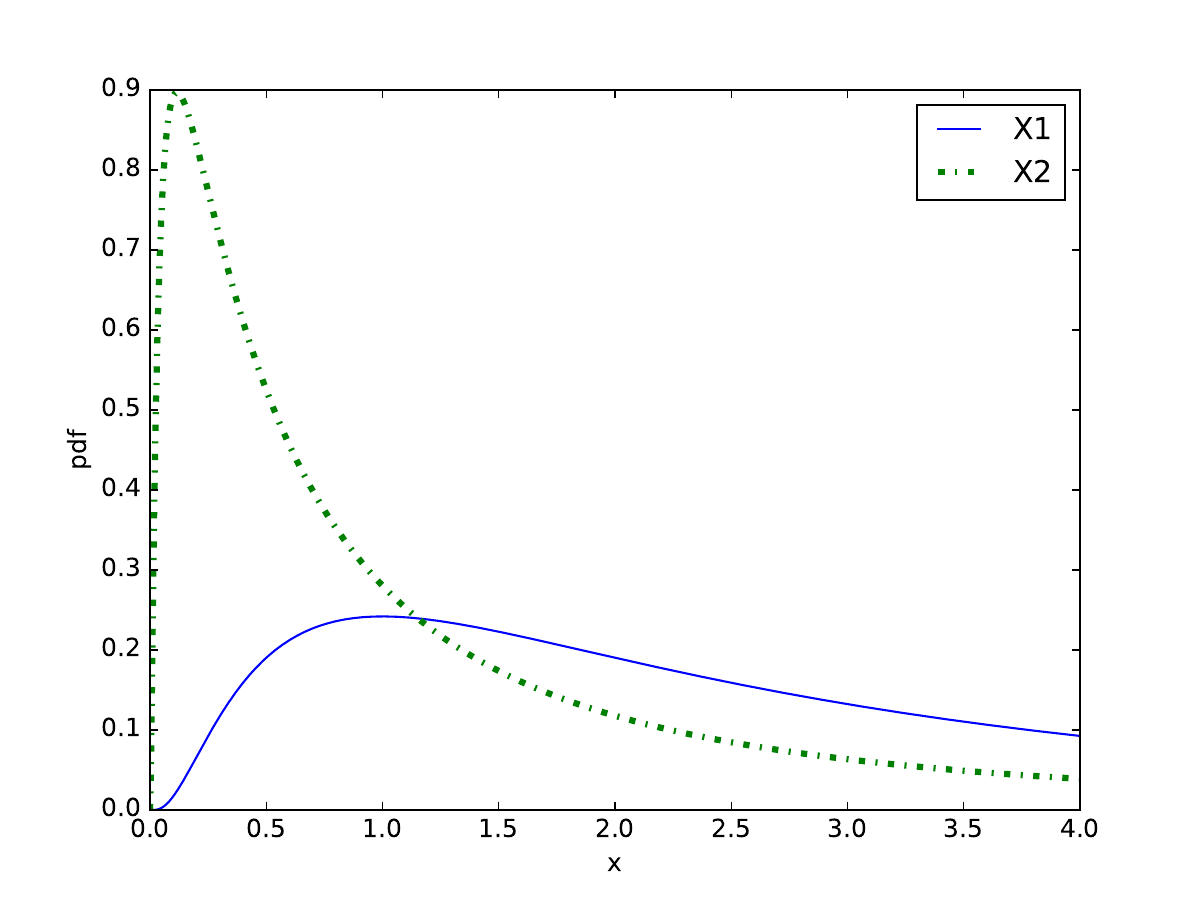}}
		\caption{Distribution of $X_1$ and $X_2$ for (a) Example \ref{ex:l1} and (b) Example \ref{ex:l2}.}
		\label{fig:ex}
	\end{figure}

	Consider objective function $c^\sigma_1$ with $\alpha = \beta = 1$. $I_1^*$ reduces to $\mathbb{E}[| X - F_X^{-1}(\frac{1}{2})|]$.  First we claim that $I_1^*(X_1) = I_1^*(X_2) = \frac{1}{4}$:
	\begin{align*}
	\mathbb{E}[| X_1 &- F_{X_1}^{-1}(\frac{1}{2})|] = \mathbb{E}[| X_1 - \frac{1}{2}|] \\
	&= \int_0^{\frac{1}{2}}(\frac{1}{2} - x)dx + \int_{\frac{1}{2}}^1(x - \frac{1}{2})dx \\
	&= \int_0^{\frac{1}{2}}(\frac{1}{2} - x)dx + \int_0^{\frac{1}{2}}xdx = \frac{1}{4}, \\
	\mathbb{E}[| X_2 &- F_{X_2}^{-1}(\frac{1}{2})|] = \mathbb{E}[| X_2 - \frac{1}{2}|] \\
	&= \int_0^{\frac{1}{2}}(\frac{1}{2} - x)4xdx + \int_{\frac{1}{2}}^1(x - \frac{1}{2})(4x - 2)dx \\
	&= \int_0^{\frac{1}{2}}(\frac{1}{2} - x)4xdx + \int_0^{\frac{1}{2}}4x^2dx = \frac{1}{4}.
	\end{align*}
	Distribution of $X_3$ can be arbitrary as long as $I_1^*(X_3) > \frac{1}{4}$ to make sure that it comes last. In order to have $I_1^*$ as the optimal index, changing the order of $X_1$ and $X_2$ should not affect the optimal value of $c^\sigma_1$. However, for the sequence $\sigma_1 \cdot \mathbf X = (X_1, X_2, X_3)$, $\inf_{\mathbf s \in \mathcal{S}}c^{\sigma_1}_1(\mathbf s) \approx 0.3946$ but sequence $\sigma_2 \cdot \mathbf X = (X_2, X_1, X_3)$ yields $\inf_{\mathbf s \in \mathcal{S}}c^{\sigma_2}_1(\mathbf s) \approx 0.3872$. 
\end{exmp}
Thus, the index $I_1^*$ is not optimal for cost function $c^\sigma_1$.

\begin{exmp}
	\label{ex:l2}
	\normalfont
	Consider objective function $c^\sigma_2$ and let $X_1 \sim \ln \mathcal{N}(1, 1)$ and $X_2 \sim \ln \mathcal{N}(\frac{1}{2}\ln (\frac{e}{e+1}), 2)$ be independent (see Figure \ref{fig:ex}).
	
	Note that $I_2^*(X_1) = I_2^*(X_2) = e^3(e - 1)$. Distribution of $X_3$ can be arbitrary as long as $I_2^*(X_3) > I_2^*(X_1) = I_2^*(X_2)$ to make sure that it comes last. In order to have $I_2^*$ as the optimal index, changing the order of $X_1$ and $X_2$ should not affect the optimal value of $c^\sigma_2$. However, for the sequence $\sigma_1 \cdot \mathbf X = (X_1, X_2, X_3)$, $\inf_{\mathbf s \in \mathcal{S}}c^{\sigma_1}_2(\mathbf s) \approx 94.158$ but sequence $\sigma_2 \cdot \mathbf X = (X_2, X_1, X_3)$ yields $\inf_{\mathbf s \in \mathcal{S}}c^{\sigma_2}_2(\mathbf s) \approx 99.096$. 
	%
	%
\end{exmp}
Thus, the index $I_2^*$ is not optimal for cost function $c^\sigma_2$.
The above leads us to the following conclusion.
\begin{thm}
	\label{thm:no opt index}
	There exists no index that yields the optimal sequence for cost functions $c^\sigma_1$ and $c^\sigma_2$.
\end{thm}
\begin{pf}
Proposition \ref{prop:optimal indexing} implies that $I_k^*$ is the only possible optimal index for cost functions $c_k, k=1,2$. But counterexamples \ref{ex:l1} and \ref{ex:l2} show that these need not be optimal. This leads us to the conclusion that optimal indices may not exist, i.e., the problem is non-indexable.~~~
\qed
\end{pf}
\begin{rem}
It is worth mentioning that Proposition \ref{prop:optimal indexing} still holds even if we restrict the space of random variables to a certain family. Therefore, although Theorem \ref{thm:no opt index} states that the sequencing problem is not indexable in general, it does not preclude the possibility of indexability in a restricted space. Nevertheless, Proposition \ref{prop:optimal indexing} ensures that one should not investigate indices other than $I_g^*$. Finding a family of distributions for which $I_g^*$ is an optimal index is still an open research problem. In particular, Example \ref{ex:l2} ensures that even if we restrict the space of random variables to exponential family, the problem remains non-indexable. In fact, we are unable to conclude about indexability if we further restrict to the exponential distribution. 
\blue{Moreover, Theorem \ref{thm:no opt index} does not exclude the possibility of existence of non-index-based optimal policies.}
\end{rem}


\subsection{Bounds on the optimal cost}

It is disappointing that contrary to long-held conjectures in the literature, the sequencing problem is non-indexable in general. Nevertheless, $I_g^*$ can be considered as a heuristic to order the random variables and achieve a suboptimal solution. We next provide lower and upper bounds on the optimum cost with $c^\sigma_1$ and $c^\sigma_2$ objective functions. Note that the Increasing order of $I_k^*, k=1, 2$ minimizes the upper bound. 
\begin{thm}
	For $k=1, 2$, the optimum cost of objective function $c^\sigma_k$ can be bounded by:
	\label{thm:upperbound}
	\begin{align}
	\label{eq:upperbound}
	\sum_{i=1}^{n-1}I_k^*(X_{\sigma(i)}) \leq \inf_{\mathbf s \in \mathcal{S}} c^\sigma_k(\mathbf s) \leq \sum_{i=1}^{n-1}(n-i)I_k^*(X_{\sigma(i)})
	\end{align}
\end{thm}
\begin{pf}
We  need two Lemmas for the proof of the theorem. Their proofs are relegated to the Appendix.
Lemma \ref{lem:tri} proves a sub-additive property of the index functions while Lemma \ref{lem:sup} is a technical lemma.
\begin{lem}\label{lem:tri}
	Let $X_1, X_2 \in \mathscr{X}^+$ be independent. Then, for $k=1, 2$
	\begin{align}
	I_k^*(X_1 + X_2) \leq I_k^*(X_1) + I_k^*(X_2).
	\end{align}
\end{lem}

\begin{lem} Assume $g(0) = 0$ and
	\label{lem:sup}
	\begin{itemize}
		\item [(i)] let $X \in \mathscr{X}^+$. Then, $\sup_{x \in \mathbb{R}}I_g^*(\max(x, X)) \leq I_g^*(X)$.
		\item [(ii)] let $X_1, X_2 \in \mathscr{X}^+$ be independent. Then, \\ $\max(I_g^*(X_1), I_g^*(X_2))\leq I_g^*(X_1 + X_2)$.
	\end{itemize}
\end{lem}

Lemmas \ref{lem:tri} and \ref{lem:sup} can now be used to bound $I_k^*(E_{\sigma(j)})$:
\begin{align}
I_k^*(E_{\sigma(j)}) &= I_k^*(\max(s_j, E_{\sigma(j-1)}) + X_{\sigma(j)}) \\
&\leq I_k^*(\max(s_j, E_{\sigma(j-1)})) + I_k^*(X_{\sigma(j)}) \\
&\leq I_k^*(E_{\sigma(j-1)}) + I_k^*(X_{\sigma(j)})
\end{align}
for $k=1, 2$ where the first and second inequality follow from Lemmas \ref{lem:tri} and \ref{lem:sup}, respectively. Using the fact that $E_{\sigma(1)} = X_{\sigma(1)}$, one can write:
\begin{align}
I_k^*(E_{\sigma(j)}) \leq \sum_{i=1}^j I_k^*(X_{\sigma(i)})
\end{align}
By lower bound in Lemma \ref{lem:sup}, $I_k^*(\max(s_j, E_{\sigma(j-1)}) + X_j) \geq I_k^*(X_j)$. Hence, $I_k^*(E_{\sigma(j)})$ can be bounded by:
\begin{align}
I_k^*(X_{\sigma(j)}) \leq I_k^*(E_{\sigma(j)}) \leq \sum_{i=1}^j I_k^*(X_{\sigma(i)}).
\end{align}

Now, to prove the upper bound let $\tilde{\mathbf s} = (\tilde{s}_2, \cdots, \tilde{s}_n)$ where $\tilde{s}_i = F_{E_{\sigma(i-1)}}^{-1}(\frac{\beta}{\alpha + \beta})$ for the case that $k=1$ and $\tilde{s}_i = \mathbb{E}[E_{\sigma(i-1)}]$ for the case that $k=2$. Note that $\tilde{s}_i$ can be calculated recursively because $E_{\sigma(i-1)}$ is a function of $\tilde{s}_2$ through $\tilde{s}_{i-1}$. We have:
\begin{align*}
\inf_\mathbf s c^\sigma_k(\mathbf s) &\leq c^\sigma_k(\tilde{\mathbf s}) = \sum_{j=2}^{n}I_k^*(E_{\sigma(j-1)}) \\
& \leq \sum_{j=1}^{n-1}\sum_{i=1}^j I_k^*(X_{\sigma(i)}) \\
& = \sum_{i=1}^{n-1}\sum_{j=i}^{n-1} I_k^*(X_{\sigma(i)}) \\
&= \sum_{i=1}^{n-1} (n - i)I_k^*(X_{\sigma(i)}).
\end{align*}
To prove the lower bound, note that $\mathbb{E}[g_k(E_{\sigma(i-1)} - s_i)] \geq I_k^*(E_{\sigma(i-1)}) \geq I_k^*(X_{\sigma(i-1)})$ where $g_1(t) = \beta(t)^+ + \alpha (-t)^+$ and $g_2(t) = t^2$. Thus,
\begin{align*}
c^\sigma_k(\mathbf s)  = \sum_{i=2}^n \mathbb{E}[g_k(E_{\sigma(i-1)} - s_i)] \geq \sum_{i=2}^{n} I_k^*(X_{\sigma(i-1)}).
\end{align*} 
\qed
\end{pf}

\begin{rem}
	Note that the upper bound and lower bound in (\ref{eq:upperbound}) coincide when $n=2$, and this is the result we already expected from Example \ref{ex:index motivation}. For general $n$, sequencing with respect to increasing order of $I_k^*$ minimizes the upper bound in (\ref{eq:upperbound}). 
\end{rem}



\section{Scheduling Problem}\label{sec:sch}

In many problems, the sequence in which to schedule is given and only the appointment times are to be determined optimally. \blue{In this section, we assume that the sequence of $n$ random variables $\mathbf X=(X_1, \cdots, X_n)$ is fixed and without loss of generality (by possibly renaming jobs) remove the notation $\sigma$ for simplicity .} We call this problem the scheduling problem. \rredd{\blue{We propose sample average approximation (SAA) as an algorithm to find the optimal appointment times and prove it is statistically consistent in the case that the objective function is convex (e.g., $c_1^\sigma$)}. This result is significant because the only assumption required for consistency of SAA is the existence of a schedule with finite cost. This assumption significantly relaxes the typical assumptions required for consistency of SAA in the literature (see e.g., Theorem 5.4 of \cite{ShDeRu09}).}

\subsection{\bblue{Existence of Solution}}
We first show that there exists a solution to the optimization problem in (\ref{eq:obj}).

\begin{thm}
	\label{thm:convexity}
	\begin{itemize}
		\item[(i)] For any particular realization of $\mathbf X$, $C(\cdot, \mathbf X)$ is nonnegative and coercive.
		\item[(ii)] $c(\cdot)$ is nonnegative, coercive and lower semi-continuous. Furthermore, if $c(\mathbf s) < \infty$ for some $\mathbf s \in \mathcal{S}$, then there exists a solution to the optimization problem in (\ref{eq:obj}) and the set of minimizers is compact.
	\end{itemize}
\end{thm}
The proof is relegated to the appendix.

One of the essential conditions in Theorem \ref{thm:convexity} is that $c(\mathbf s) < \infty$ for some $\mathbf s \in \mathcal{S}$. The question is how to check whether this condition is satisfied. Should we explore the entire set $\mathcal{S}$ in the hope of finding such $\mathbf s$? Let's illuminate this condition: First of all it is easy to see that for $p \geq 1$ and $g(t) = |t|^p$, this condition is equivalent to $X_i \in L^p$ (i.e., $\mathbb{E}[|X_i|^p] < \infty$) for $i=1, \cdots, n-1$. This is also true for some other variations where $g$ is a piecewise function of the form $|\cdot|^p$ such as $g_1$ and $g_3$ in Examples \ref{ex:c1} and \ref{ex:g3}. Moreover, if $c(\mathbf s) < \infty$ for some $\mathbf s \in \mathcal{S}$, it is finite for all $\mathbf s \in \mathcal{S}$. It is mainly due to the fact that $L^p$ is a vector space. Therefore, in such cases, there is no need to explore the set $\mathcal{S}$. However, for general $g$, the set $\{X \in \mathscr{X}^+ \mid \mathbb{E}[g(X)]< \infty]\}$ may not be a vector space (see Birnbaum-Orlicz space, \cite{BiOr31})  and $c(\mathbf s)$ may be infinite for some  $\mathbf s$. In that case, random exploration may yield $\mathbf s \in \mathcal{S}$ such that $c(\mathbf s) < \infty$.

\subsection{\bblue{Sample Average Approximation}}
The next question is how to calculate the optimal appointment times. Theorem \ref{thm:convexity} assures that there exists an optimal schedule under mild condition. However, calculating expectation is very costly in our problem \bblue{due to the convolution nature of the distribution of the service completion times. In fact, for a given schedule $\mathbf s$, distribution of $E_{\sigma(i)}$ is convolution of distributions of $\max (s_i, E_{\sigma(i-1)})$ and $X_i$}. An alternative is to use sample average approximation (SAA) to approximate the optimization problem. \bblue{SAA is a well studied topic in stochastic programming (see for example, \cite{ShDeRu09,royset2013sample})}. \bblue{In the following, we discuss SAA and provide a theoretical guarantee for convergence of the solution in stochastic appointment scheduling problem.} \blue{We assume that
\begin{assum}
\label{assum: convexity of C}
For any realization of $\mathbf X$, $C(\cdot, \mathbf X)$ is convex.
\end{assum}
This assumption holds for the $l_1$-type objective function (see Proposition \ref{prop: convexity of g1} in Appendix A) which is widely considered in the literature. However, it does not hold for the case of $l_2$-type objective (see Example \ref{ex: non-convexity of g2} in Appendix A).
}

Let $(\mathbf X^j)_{j=1}^m$ be an \blue{independently and identically distributed} (i.i.d.) random sample of size $m$ for durations $\mathbf X$ and define
\begin{align}
C_m(\mathbf s) = \frac{1}{m}\sum_{j=1}^m C(\mathbf s, \mathbf X^j)
\end{align}
Instead of solving the optimization problem in Equation \ref{eq:obj}, we're going to solve
\begin{align}
\label{eq:saa}
\inf_{\mathbf s \in \mathcal{S}} C_m(\mathbf s).
\end{align}
Convexity and coercivity of $C(\cdot, \mathbf X)$  implies convexity and coercivity of $C_m(\cdot)$. Therefore, there exists a solution to the optimization problem in (\ref{eq:saa}). In addition, Strong Law of Large Numbers implies that for each $\mathbf s$, $C_m(\mathbf s) \to c(\mathbf s)$ a.s. as $m \to \infty$. Nevertheless, optimization over the set $\mathcal{S}$ requires some stronger result to guarantee $\inf_{\mathbf s \in \mathcal{S}} C_m(\mathbf s) \to \inf_{\mathbf s \in \mathcal{S}}c(\mathbf s)$ a.s. as $m \to \infty$. Moreover, it would be useful to see if the set of minimizers of the SAA also converges to the set of true minimizers in some sense. To reach that goal, we need the following definition of deviation for sets (see equation (7.4)  in \cite{ShDeRu09}).
\begin{defn}
	Let $(M, d)$ be a metric space and $A, B \subseteq M$. We define distance of $a \in A$ from $B$ by
	\begin{align}
	\text{dist}(a, B) := \inf \{d(a, b) \mid b \in B\}
	\end{align}
	and deviation of $A$ from $B$ by
	\begin{align}
	\mathbb{D}(A, B) := \sup_{a \in A} \text{dist}(a, B).
	\end{align}
\end{defn}
Note that $\mathbb{D}(A, B) = 0$ implies $A \subseteq \text{cl}(B)$ (i.e. $A$ is a subset of closure of $B$ with respect to $M$). The next theorem guarantees that SAA is a consistent estimator for the scheduling problem.
\begin{thm}
	\label{thm:saa}
	Suppose Assumption \ref{assum: convexity of C} holds and $c(\mathbf s) < \infty$ for some $\mathbf s \in \mathcal{S}$ and let $S^* = \emph{arginf}_{\mathbf s \in \mathcal{S}}c(\mathbf s)$ and $S_m^* = \emph{arginf}_{\mathbf s \in \mathcal
		S}C_m(\mathbf s)$. Then,
	$\inf_{\mathbf s \in \mathcal{S}}C_m(\mathbf s) \to \inf_{\mathbf s \in \mathcal{S}}c(\mathbf s)$ and $\mathbb{D}(S_m^*, S^*) \to 0$ a.s. as $m \to \infty$.
\end{thm}
The proof is available in the appendix.

Theorem \ref{thm:saa} proves the consistent behavior of SAA as the number of samples tends to infinity. Let's now observe how it behaves in terms of bias. For any $\mathbf s' \in \mathcal{S}$, we can write $\inf_{\mathbf s \in \mathcal{S}} C_m(\mathbf s) \leq C_m(\mathbf s')$. By taking expectation and then minimizing over $\mathbf s'$, we conclude that $\mathbb{E}[\inf_{\mathbf s \in \mathcal{S}} C_m(\mathbf s)] \leq \inf_{\mathbf s \in \mathcal{S}}\mathbb{E}[C_m(\mathbf s)]$. Since samples are i.i.d., $\mathbb{E}[C_m(\mathbf s)] = c(\mathbf s)$. Therefore, $\mathbb{E}[\inf_{\mathbf s \in \mathcal{S}} C_m(\mathbf s)] \leq \inf_{\mathbf s \in \mathcal{S}}c(\mathbf s)$ which means SAA is negatively biased. Does this bias decrease as the number of samples increases? The answer is affirmative. \bblue{Theorem 2 in \cite{MaMoWo99} proves that $\mathbb{E}[\inf_{\mathbf s \in \mathcal{S}} C_m(\mathbf s)] \leq \mathbb{E}[\inf_{\mathbf s \in \mathcal{S}} C_{m+1}(\mathbf s)]$.}

\section{Numerical Results}
\label{sec:numerical}

\bblue{It has become a standard practice to evaluate performance on operating room data due to the immediate application of stochastic appointment scheduling in healthcare. \cite{DeViVo07} used real surgery scheduling data collected at Fletcher Allen Health Care of New York. In this paper, we consider surgery scheduling dataset from Keck hospital of USC.}

\bblue{The dataset includes 38,000 surgeries performed in 25 operating rooms over the course of 3 years. More than 800 different procedure types performed by 200 surgeons. \rredd{Surgeries with the same procedure type performed by the same surgeon are assumed to be samples of the same distribution. Our numerical analysis is restricted to those distributions that have at least 30 samples.}
We stick to 30 samples because we observed that they are sufficient for a close enough SAA of the optimal solution. This is much fewer than the theoretically required number of samples given by \cite{BeLeQu12}.}
\rredd{In some practical scenarios, there are not enough samples to directly apply SAA. In such scenarios, similar cases based on the nature of the procedure type can be aggregated to build distributions with enough number of samples. In this paper, we focus on the surgeon-procedure pairs that have enough number of samples.}

\begin{figure}
	
	\centering
	\includegraphics[scale=.52]{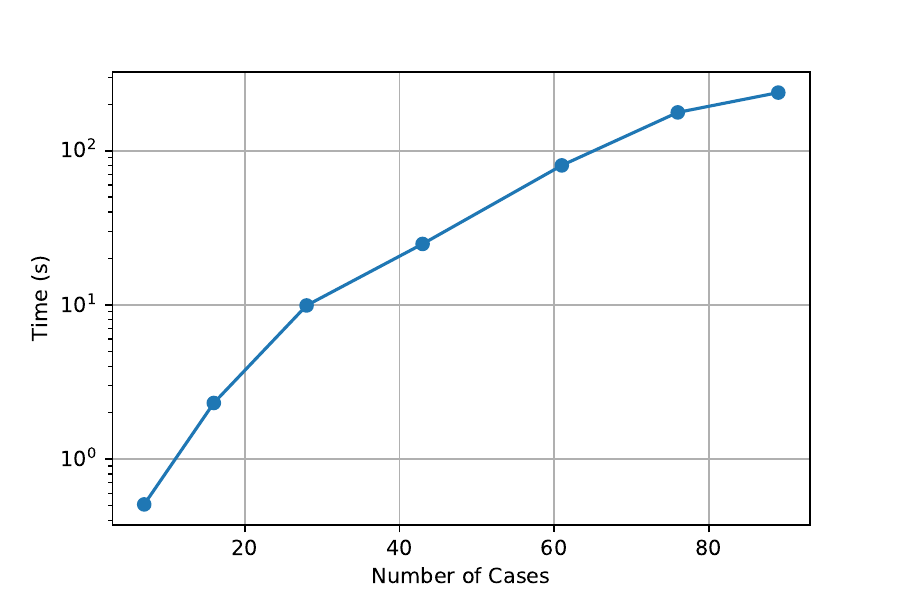}	
	\caption{SAA running time (in seconds) to find approximate optimal schedule for a given sequence. 30 samples/job are used for SAA though no appreciable difference even if 10x more samples used.}
	\label{fig:time}	
\end{figure}

\bblue{We first show that given a sequence, SAA-based optimization algorithm is fast enough for all practical purposes to find an approximate solution. To do so, we use the Powell method (\cite{brent2013algorithms}) to solve the SAA-based optimization problem numerically. The experiments are performed in Python on a 2015 Macbook Pro with 2.7 GHz Intel Core i5 processor and 16 GB 1867 MHz DDR3 memory. Figure \ref{fig:time} confirms that appointments for a given sequence of $n = 80$ jobs can be calculated in about 3 minutes. 
Moreover, we observed that changing the number of samples from 10 to 300 does not change the run time of the SAA-based optimization significantly.}

\begin{figure}
	\centering     
	\subfigure[]{\label{fig:boundsa}\includegraphics[width=80mm]{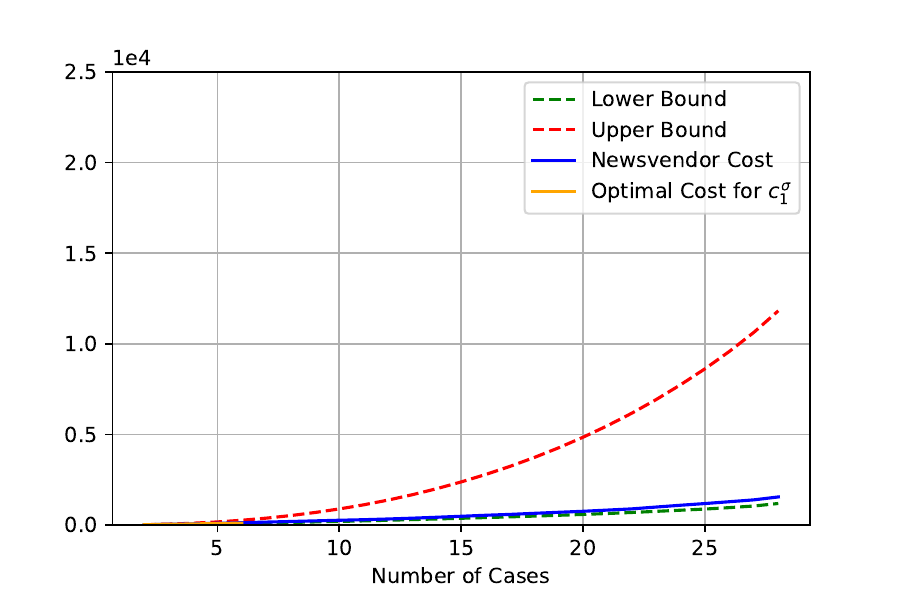}}
	\subfigure[]{\label{fig:boundsb}\includegraphics[width=80mm]{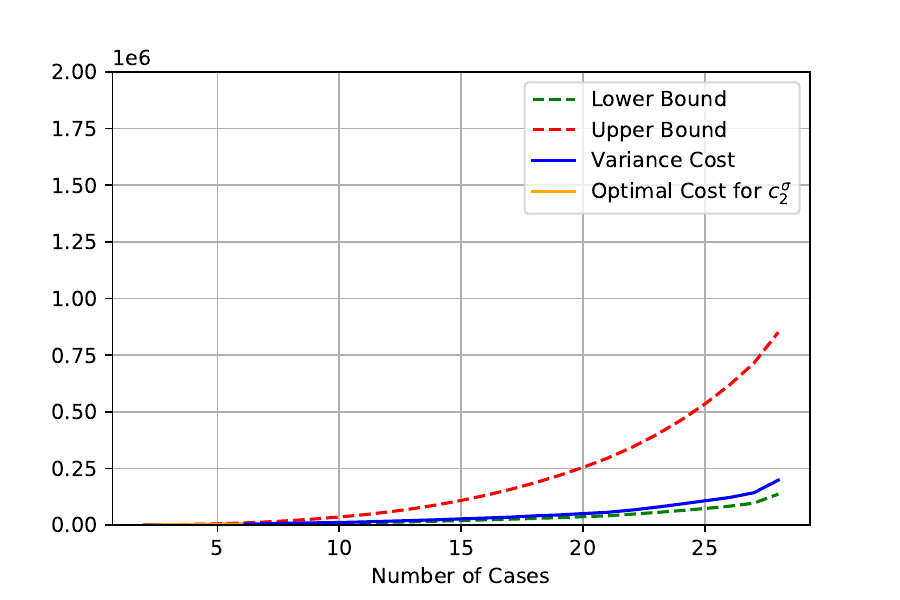}}
	\caption{Upper and lower bounds on optimal cost for (a) $c^\sigma_1$ and (b) $c^\sigma_2$ cost functions. As shown in the figure, upper bound is quite loose on USC Keck dataset. Table \ref{tab: zoom in} provides numerical values for $n \leq 6$ to compare optimal sequence with index-based heuristic policy.}
	\label{fig:bounds}
\end{figure}

Secondly, the bounds provided in Theorem \ref{thm:upperbound} are evaluated. Bounds in Theorem \ref{thm:upperbound} are for general distribution and may be useful in the worst case scenarios. However, Figure \ref{fig:bounds} shows that the upper bound is loose as the number of jobs increases on Keck dataset. \blue{The upper bound of Theorem \ref{thm:upperbound} uses the complete delay propagation through the schedule, i.e., potential cost of each job affects all the future jobs equally. Although this situation might arise in the worst case, we've observed that on Keck dataset, it does not happen. Indeed, the gaps between jobs prevents the delay to have full effect on subsequent jobs.} 

\bblue{Non-indexability shown in Theorem \ref{thm:no opt index} is for general distribution. One might wonder if non-indexability is actually observed in practice. We verify that the optimal sequence is indeed different from the one given by heuristic policies (\blue{see Table \ref{tab: zoom in}}) using Keck dataset. Newsvendor and LVF indexes are considered as heuristic policies for $c^\sigma_1$ and $c^\sigma_2$ objective functions, respectively since they are the only possible candidates to return the optimal sequence (Proposition \ref{prop:optimal indexing}). The true optimal sequence is calculated by comparing all $n!$ choices. In operating room scheduling, the number of surgeries performed in a typical day hardly exceeds 6 which leaves the door open for exhaustive search to find the optimal sequence. However, other applications such as outpatient clinics have much larger number of jobs and it may not be feasible to exhaustively search over all possible sequences.}

\begin{figure}
	\centering     
	\includegraphics[width=80mm]{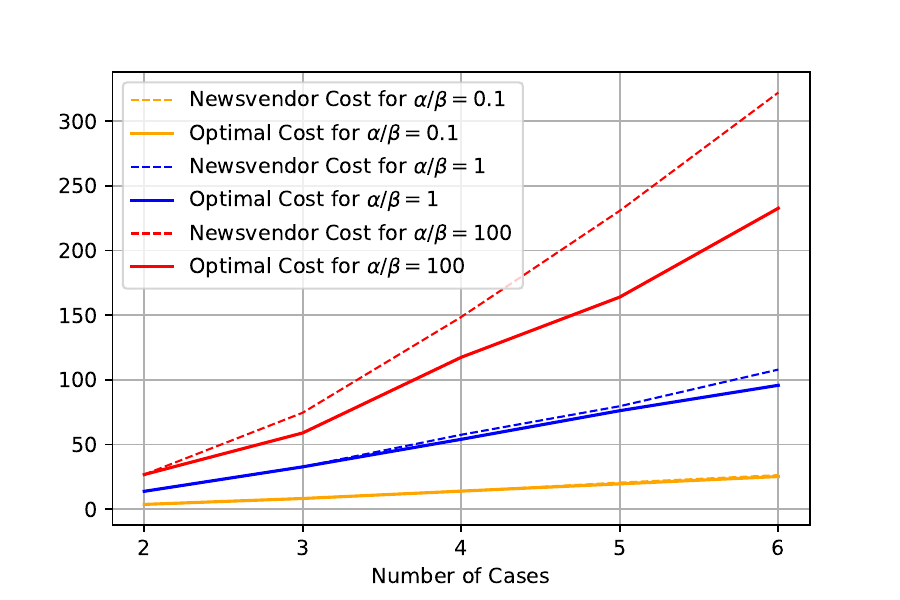}
	\caption{\blue{Optimality gap of newsvendor index increases as the ratio of $\alpha/\beta$ increases on Keck dataset. The optimal sequence is found by exhaustive search over all $n!$ possible sequences. $\beta = 1$ is fixed and $\alpha$ changes from $0.1$ to $100$. The dashed lines show the cost for the sequence obtained by least newsvendor first. The cost of the optimal sequence is shown by solid lines.}}
	\label{fig:sensitivity}
\end{figure}

\blue{Cost function $c^\sigma_1$ depends on the idle time and delay per unit costs $\alpha$ and $\beta$. \cite{MaMaZiKr18} analyzed how newsvendor index outperforms variance in different regimes of these parameters. In Figure \ref{fig:sensitivity}, we evaluate the gap between newsvendor index and the optimal sequence as $\alpha$ and $\beta$ change. The optimal sequence is obtained by exhaustive search over all $n!$ possible sequences. It can be seen that as the ratio of $\alpha/\beta$ increases, the sub-optimality gap of newsvendor index increases on Keck dataset.}

\begin{table*}
\caption{\blue{Non-optimality of least newsvendor first for $c^\sigma_1$ and least variance first for $c^\sigma_2$. Optimal sequence found by exhaustive search is different from the sequence given by heuristic index-based policies.}}
\label{tab: zoom in}
\begin{center}
\renewcommand{\arraystretch}{1.3}
\begin{tabular}{ |l|c|c|c|c|c|c| }
\hline
&& \multicolumn{5}{ |c| }{$n$ (Number of jobs)} \\
\cline{3-7}
 && 2 & 3 & 4 & 5 & 6\\ \hline
\multirow{4}{*}{$c^\sigma_1$} & \textbf{lower bound} & 14.9 & 32.8 & 51.2 & 71.4 & 92.3\\ \cline{2-7}
& \textbf{optimal cost} & 14.9 & 36.6 & 51.6 & 79.0 & 105.3 \\ \cline{2-7}
& \textbf{newsvendor cost} & 14.9 & 36.6 & 64.4 & 95.4 & 126.5 \\ \cline{2-7}
& \textbf{upper bound} & 14.9 & 47.7 & 98.9 & 170.4 & 262.7 \\ \cline{2-7}
 \hline \hline
\multirow{4}{*}{$c^\sigma_2$} & \textbf{lower bound} & 368.0 & 817.5 & 1534.3 & 2430.0 & 3562.3 \\ \cline{2-7}
& \textbf{optimal cost} & 368.0 & 1036.3 & 1763.8 & 2760.1 & 3912.5 \\ \cline{2-7}
& \textbf{variance cost} & 368.0 & 1081.7 & 1923.1 & 2853.7 & 3939.1 \\ \cline{2-7}
& \textbf{upper bound} & 368.0 & 1185.5 & 2719.9 & 5149.8 & 8712.2 \\ \hline
\end{tabular}
\end{center}
\end{table*}

\section{Conclusions}
\label{sec:conclusion}

In this paper, we considered the optimal stochastic appointment scheduling problem. Each job potentially has a different service time distribution and the objective is to minimize the expectation of a function of idle time and start-time delay. There are two sub-problems. (i) The sequencing problem: the optimal sequence in which to schedule the jobs. We show that this problem in general is non-indexable. (ii) The scheduling problem: finding the optimal appointment times given a sequence or order of jobs.  We show that there exists a solution to the scheduling problem. Moreover, the $l_1$-type objective function is convex. Further, we give a sample average approximation-based algorithm that yield an approximately optimal solution which is asymptotically consistent. 

It has been an open problem for many years to find the index that yields the optimal sequence of jobs. Following the work of \cite{We90}, who showed that Least Variance First (LVF) is optimal for two cases for specific distributions, it had been conjectured that the problem is indexable and LVF may be optimal for the general problem with the $l_1$-type objective. In fact, several simulation studies and approximation algorithms are based on such policies. In this paper, we have settled the open question of the optimal index-type policy, namely that the problem is non-indexable in general, and no such index exists. Indeed, we show that if the problem is indexable, then a `Newsvendor index' would be optimal for the $l_1$ cost objective, a variance index would be optimal for $l_2$ objective, and we also give form of an index $I^*_g$ that would be optimal for a generalized cost function $g$. But we provide counterexamples that show that an optimal index-based policy does not  exist for some problems. It is quite possible that the problem is indexable for specific distribution classes. That remains an open research question.
 
\bibliographystyle{unsrt}
\bibliography{sch-survey}

\appendix
\section{Proofs}
\label{sec: appendixa}
\textbf{Proof of Lemma \ref{lem:tri}:}
For $k=2$ the statement is obvious. For $k=1$, using the fact that $(a + b)^+ \leq a^+ + b^+$ for $a, b \in \mathbb{R}$, we have:
\begin{align*}
&I_1^*(X_1 + X_2) \\
&= \inf_s \mathbb{E}[\alpha(s - X_1 - X_2)^+ + \beta (X_1 + X_2 - s)^+] \\
& \leq \mathbb{E}[\alpha(F_{X_1}^{-1}(\frac{\beta}{\alpha + \beta}) + F_{X_2}^{-1}(\frac{\beta}{\alpha + \beta}) - X_1 - X_2)^+ \\
&+ \beta (X_1 + X_2 - F_{X_1}^{-1}(\frac{\beta}{\alpha + \beta}) - F_{X_2}^{-1}(\frac{\beta}{\alpha + \beta}))^+] \\
& \leq \mathbb{E}[\alpha(F_{X_1}^{-1}(\frac{\beta}{\alpha + \beta}) - X_1)^+ + \beta (X_1 - F_{X_1}^{-1}(\frac{\beta}{\alpha + \beta}))^+] \\
&+ \mathbb{E}[\alpha(F_{X_2}^{-1}(\frac{\beta}{\alpha + \beta}) - X_2)^+ + \beta (X_2 - F_{X_2}^{-1}(\frac{\beta}{\alpha + \beta}))^+] \\
&= I_1^*(X_1) + I_1^*(X_2) ~~~~~
\end{align*}
\qed

\textbf{Proof of Lemma \ref{lem:sup}:}
(i) We can write $g(t) = g_r(t) + g_l(t)$ where $g_r(t) = g(t^+)$ and $g_l(t) = g(-(-t)^+)$ capture $g$ for positive and negative values of $t$, respectively. Since $g$ is nonnegative, convex and $g(0) = 0$, we can conclude that $g_r$ is nondecreasing and $g_l$ is nonincreasing. Moreover, $I_g^*(X) = \inf_{s \geq 0} \mathbb{E}[g(X-s)] = \inf_{s \in \mathbb{R}} \mathbb{E}[g(X-s)]$. Suppose $s^*$ is a minimizer for $I_g^*$ and let $\mathcal{X} = \{x \in \mathbb{R}: x \leq s^*\}$. We prove the lemma for $x \in \mathcal{X}$ and $x \notin \mathcal{X}$ separately.

Let $x \in \mathcal{X}$. We can write:
\begin{align*}
&I_g^*(\max (x, X)) = \inf_s \mathbb{E}[g(\max (x, X) - s)] \\
&= \inf_s \mathbb{E}[g_r(\max (x, X) - s) + g_l(\max (x, X) - s)] \\
& \leq \inf_{s \geq x} \mathbb{E}[g_r(\max (x, X) - s) + g_l(\max (x, X) - s)] \\
& = \inf_{s \geq x} \mathbb{E}[g_r(X - s) + g_l(\max (x, X) - s)] \\
& \leq \inf_{s \geq x} \mathbb{E}[g_r(X - s) + g_l(X - s)] \\
&= \inf_{s \geq x}\mathbb{E}[g(X - s)] \\
&= \mathbb{E}[g(X - s^*)] = I_g^*(X)
\end{align*}
For the case that $x \notin \mathcal{X}$, we can write:
\begin{align*}
&I_g^*(\max (x, X)) = \inf_s \mathbb{E}[g(\max (x, X) - s)] \\
&= \inf_s \mathbb{E}[g_r(\max (x, X) - s) + g_l(\max (x, X) - s)] \\
& \leq \inf_{s < x} \mathbb{E}[g_r(\max (x, X) - s) + g_l(\max (x, X) - s)] \\
& = \inf_{s < x} \mathbb{E}[g_r(\max (x, X) - s)] \\
&= \mathbb{E}[g_r(\max (x, X) - x)] \\
&= \mathbb{E}[g_r(X - x)] \\
&\leq \mathbb{E}[g_r(X - s^*)] \\
&\leq \mathbb{E}[g_r(X - s^*) + g_l(X - s^*)] \\
&= \mathbb{E}[g(X - s^*)] =  I_g^*(X)
\end{align*}
(ii) Note that $I_g^*(X) = \inf_{s \geq 0} \mathbb{E}[g(X-s)] = \inf_{s \in \mathbb{R}} \mathbb{E}[g(X-s)]$. To prove $\max(I_g^*(X_1), I_g^*(X_2))\leq I_g^*(X_1 + X_2)$, by symmetry, suffices to prove $I_g^*(X_1)\leq I_g^*(X_1 + X_2)$. Let $x_2 \geq 0$, we have:
\begin{align*}
I_g^*(X_1) &= \inf_s \mathbb{E}[g(X_1 - s)] \\
&= \inf_s \mathbb{E}[g(X_1 + x_2 - s)] \\
&= \inf_s \mathbb{E}[g(X_1 + X_2 - s) \mid X_2 = x_2] \\
&= \inf_s \phi(s, x_2)
\end{align*}
where $\phi(s, x_2) = \mathbb{E}[g(X_1 + X_2 - s \mid X_2 = x_2)]$. The above equality holds for any value of $x_2 \geq 0$. Hence, $I_g^*(X_1) \leq \phi(s, X_2)$ for any $s \in \mathbb{R}$. Therefore, $I_g^*(X_1) \leq \mathbb{E}[\phi(s, X_2)] = \mathbb{E}[g(X_1 + X_2 - s)]$ by smoothing property of conditional expectation. Thus, $I_g^*(X_1) \leq \inf_s \mathbb{E}[g(X_1 + X_2 - s)] = I_g^*(X_1 + X_2)$.
\qed

\textbf{Proof of Theorem \ref{thm:convexity}.}
(i) Since $g$ is nonnegative, it is obvious that $C(\cdot,\mathbf X)$ is also nonnegative.

To prove coercivity of $C(\cdot, \mathbf X)$, let $(\mathbf s^m)_{m \geq 1} \subseteq \mathbb{R}^{n-1}$ be a sequence such that $\|\mathbf s^m\| \to \infty$. We need to show that $C(\mathbf s^m, \mathbf X) \to \infty$ as $m \to \infty$. Let $j$ be the smallest integer such that $\|s_j^m\| \to \infty$. Note that for any particular realization of $\mathbf X$, there exists $M \in \mathbb{R}$ such that $|E_{j-1}^m| \leq M$ for all $m$ where $E_{j-1}^m$ denotes finish time of job $j-1$ with schedule $\mathbf s^m$. By triangle inequality,
\begin{align}
|E_{j-1}^m - s_j^m| \geq |s_j^m| - |E_{j-1}^m| \geq |s_j^m| - M \to \infty
\end{align}
as $m \to \infty$. Coercivity of $g$ implies that $g(E_{j-1}^m - s_j^m) \to \infty$ as $m \to \infty$. On the other hand, since $g$ is nonnegative we can write $C(\mathbf s^m, \mathbf X) \geq g(E_{j-1}^m - s_j^m)$ for all $m$. Hence, $C(\mathbf s^m, \mathbf X) \to \infty$ as $m \to \infty$.

(ii) Clearly, $c(\cdot)$ is nonnegative. To prove coercivity, let $(\mathbf s^m)_{m \geq 1}$ be as defined in the previous part, by Fatou's Lemma and coercivity of $C(\cdot, \mathbf X)$ we have:
\begin{align*}
\liminf_m c(\mathbf s^m) &= \liminf_m \mathbb{E}[C(\mathbf s^m, \mathbf X)] \\
&\geq \mathbb{E}[\liminf_m C(\mathbf s^m, \mathbf X)] = \infty.
\end{align*}
To prove lower semi-continuity, let $(\mathbf s^k)_{k \geq 1} \subseteq \mathbb{R}^{n-1}$ be a sequence converging to $\mathbf s \in \mathbb{R}^{n-1}$. By Fatou's Lemma we can write
\begin{align*}
\liminf_k c(\mathbf s^k) &= \liminf_k \mathbb{E}[C(\mathbf s^k, \mathbf X)]  \\
&\geq \mathbb{E}[\liminf_k C(\mathbf s^k, \mathbf X)] \geq \mathbb{E}[C(\mathbf s,\mathbf X)] = c(\mathbf s).
\end{align*}

\blue{Since $c$ is coercive and $c(\mathbf s) < \infty$ for some $\mathbf s \in \mathcal{S}$, without loss of generality we can assume that the minimization is over a compact set. Moreover, $c(\mathbf s)$ is lower semi-continuous. Thus, the set of minimizers is nonempty and compact.}
\qed

\textbf{Proof of Theorem \ref{thm:saa}.}
Define the extended real valued functions
\begin{align*}
\bar{C}_m(\mathbf s) &= C_m(\mathbf s) + \mathbb{I}_\mathcal{S}(\mathbf s) \\
\bar{c}(\mathbf s) &= c(\mathbf s) + \mathbb{I}_\mathcal{S}(\mathbf s)
\end{align*}
where
\begin{align*}
\mathbb{I}_\mathcal{S}(\mathbf s) = \begin{cases}
0, & \text{if } \mathbf s \in \mathcal{S}  \\
+\infty, & \text{Otherwise} 
\end{cases}
\end{align*}
Note that $\bar{C}_m, \bar{c}$ are nonnegative, convex and lower semicontinuous because $C_m, c$ are lower semicontinuous and $\mathcal{S}$ is closed and convex. By Theorem 2.3 of \cite{ArWe94} \bblue{(see Appendix B)}, $\bar{C}_m(\cdot)$ epi-converges to $\bar{c}(\cdot)$ (denoted by $\bar{C}_m(\cdot) \xrightarrow{e} \bar{c}(\cdot)$) for a.e. $\omega \in \Omega$.

Note that $S^* = \text{arginf}_{\mathbf s \in \mathcal{S}}c(\mathbf s) =  \text{arginf}_{\mathbf s \in \mathbb{R}^{n-1}}\bar{c}(\mathbf s)$ and $S_m^* = \text{arginf}_{\mathbf s \in \mathcal{S}}C_m(\mathbf s) = \text{arginf}_{\mathbf s \in \mathbb{R}^{n-1}}\bar{C}_m(\mathbf s)$. Since $c(\mathbf s) < \infty$ for some $\mathbf s \in \mathcal{S}$, by Theorem \ref{thm:convexity} we know that $S^*$ is nonempty and compact. Let $K$ be a compact subset of $\mathbb{R}^{n-1}$ such that $S^*$ lies in the interior of $K$. Let $\hat{S}_m^* = \text{arginf}_{\mathbf s \in K}\bar{C}_m(\mathbf s)$. We first show that for a.e. $\omega \in \Omega$, $\hat{S}_m^*$ is nonempty for large enough $m$. Let $\mathbf s^* \in S^*$ and consider $\omega \in \Omega$ for which $\bar{C}_m(\cdot) \xrightarrow{e} \bar{c}(\cdot)$. By definition of epi-convergence, $\limsup_m \bar{C}_m(\mathbf s_m) \leq \bar{c}(\mathbf s^*)$ for some $\mathbf s_m \to\mathbf  s^*$. Therefore, there exists $M \geq 1$ such that for $m \geq M$, $\bar{C}_m(\mathbf s_m) \leq \bar{c}(\mathbf s^*) + 1 < \infty$. Moreover, it follows from $\mathbf s_m \to \mathbf s^*$ that for large enough $m$, $\mathbf s_m$ lies in the interior of $K$. Since $\bar{C}_m(\cdot)$ is convex and lower semicontinuous and $K$ is compact, $\hat{S}_m^*$ is nonempty a.s. (see Appendix B for Proposition 2.3.2 of \cite{BeNeOz03}).

Now, let us show that $\mathbb{D}(\hat{S}_m^*, S^*) \to 0$ a.s. Consider $\omega \in \Omega$ for which $\bar{C}_m(\cdot) \xrightarrow{e} \bar{c}(\cdot)$. We claim that for such $\omega$, $\mathbb{D}(\hat{S}_m^*, S^*) \to 0$. Assume by contradiction that $\mathbb{D}(\hat{S}_m^*, S^*) \not\to 0$. Thus, there exists $\epsilon > 0$ and $\mathbf y_m \in \hat{S}_m^*$ (for large enough $m$) such that $\text{dist}(\mathbf y_m, S^*) \geq \epsilon$. Let $\mathbf y_{m_l} \to \mathbf y$ be a convergent subsequence of $(\mathbf y_m)_{m \geq 1}$. Such a subsequence exists because $K$ is compact. It follows from $\text{dist}(\mathbf y_m, S^*) \geq \epsilon$ that $\mathbf y \notin S^*$. On the other hand, Proposition 7.26 of \cite{ShDeRu09} \bblue{(see Appendix B)} implies that $\mathbf y \in \text{arginf}_{\mathbf s \in K}\bar{c}(\mathbf s) = S^*$ which is a contradiction.

Note that $S^*$ is in the interior of $K$. It follows from $\mathbb{D}(\hat{S}_m^*, S^*) \to 0$ that for large enough $m$, $\hat{S}_m^*$ lies in the interior of $K$. Hence, $\hat{S}_m^*$ is a local minimizer. Convexity of $\bar{C}_m(\cdot)$ implies that $\hat{S}_m^*$ is a global minimizer i.e. $\hat{S}_m^* = S_m^*$. Therefore, $\mathbb{D}(S_m^*, S^*) \to 0$ a.s. as $m \to \infty$.

It remains to prove that $\inf_{\mathbf s \in \mathcal{S}}C_m(\mathbf s) \to \inf_{\mathbf s \in \mathcal{S}}c(\mathbf s)$ a.s. Fix $\omega \in \Omega$ for which $\bar{C}_m(\cdot) \xrightarrow{e} \bar{c}(\cdot)$ and let $\mathbf s_m^* \in S_m^*$ be a convergent sequence. Such a sequence exists because for large enough $m$, $S_m^*$ falls inside the compact set $K$. Then, by Proposition 7.26 of \cite{ShDeRu09} (see below), $\inf \bar{C}_m(\mathbf s) \to \inf \bar{c}(\mathbf s)$ or equivalently, $\inf_{\mathbf s \in \mathcal{S}}C_m(\mathbf s) \to \inf_{\mathbf s \in \mathcal{S}}c(\mathbf s)$. 
\qed

\blue{
\begin{prop}
\label{prop: convexity of g1} Let $\mathbf X$ be a fixed sequence of jobs and $C_1(\cdot, \mathbf X) = \sum_{i=2}^n \left[\alpha(s_i - E_{i-1})^+ \\ + \beta(E_{i-1} - s_i)^+\right]$ as in Example \ref{ex:c1}. For any realization of $\mathbf X$, $C_1(\cdot, \mathbf X)$ is convex and thus, $c_1(\cdot) = \mathbb{E}[C_1(\cdot, \mathbf X)]$ is also convex.
\end{prop}
}
\begin{pf}
\blue{
We proceed by writing $C_1(\cdot, \mathbf X)$ as maximum of $2^{n-1}$ affine functions. Since an affine function is convex, so is the maximum. To define these functions, we first split $\mathbb{R}^{n-1}$ into $2^{n-1}$ regions and then define an affine function in each region. These functions are then extended to the entire $\mathbb{R}^{n-1}$. The detail is given in the following.
}

\blue{
Fix a realization of $\mathbf X$ and note that for any schedule $\mathbf s = (s_2, \cdots, s_n) \in \mathbb{R}^{n-1}$, either $(s_i - E_{i-1})^+ > 0$ or $(E_{i-1} - s_i)^+ > 0$ for $i = 2, \cdots, n$. Let $b_i$ be a binary variable that indicates which of the two happens. More specifically, $b_i = 1$ if $s_i \geq E_{i-1}$, and $b_i = 0$ if $s_i < E_{i-1}$. These binary variables are used to split $\mathbb{R}^{n-1}$ into $2^{n-1}$ regions $R_{b_2, \cdots, b_n}$. More precisely, if $b_i = 1$, then $s_i \geq E_{i-1}$ denotes the range of $s_i$ in $R_{b_2, \cdots, b_n}$ and if $b_i=0$, then $s_i < E_{i-1}$ determines its range. For example for $n=4$, region $R_{101}$ would be
\begin{align*}
R_{101} := \{(s_2, s_3, s_4) \in \mathbb{R}^3 \mid s_2 \geq E_1, s_3< E_2, s_4 \geq E_3\}.
\end{align*}
Corresponding to each region, one can define a function $\bar f_{b_2, \cdots, b_n}: R_{b_2, \cdots, b_n} \to \mathbb{R}$ that consists of sum of $n-1$ terms associated with each $b_i$. If $b_i = 1$, then the corresponding term would be $\alpha (s_i - E_{i-1})$ and if $b_i = 0$, it would be $\beta (E_{i-1} - s_i)$. For example for $n=4$, $\bar f_{101}$ would be
\begin{align}
\label{eq: tmp fbar}
&\bar f_{101}(\mathbf s) := \alpha(s_2 - E_{1}) + \beta (E_2 - s_3) + \alpha (s_4 - E_3)  \nonumber \\
&= \alpha(s_2 - X_1) + \beta(s_2 + X_2 - s_3) + \alpha (s_4 - s_2 - X_2).
\end{align}
Note that $C_1(\mathbf s, \mathbf X) =  \bar f_{b_2, \cdots, b_n}(\mathbf s)$ on $R_{b_2, \cdots, b_n}$. Moreover, restricting the domain of $\bar f_{b_2, \cdots, b_n}$ to $R_{b_2, \cdots, b_n}$ allowed us to write the last equality in \eqref{eq: tmp fbar} which can now be used for an affine extension to the entire $\mathbb{R}^{n-1}$. Let $f_{b_2, \cdots, b_n}: \mathbb{R}^{n-1} \to \mathbb{R}$ be such an extension. We claim that $C_1(\mathbf s, \mathbf X) = \max_{b_2, \cdots, b_n} f_{b_2, \cdots, b_n}(\mathbf s)$ for all $\mathbf s \in \mathbb{R}^{n-1}$ and thus convex. To prove this claim, it suffices to show that on $R_{b_2, \cdots, b_n}$, $\bar f_{b_2, \cdots, b_n}(\mathbf s) \geq f_{b'_2, \cdots, b'_n}(\mathbf s)$ (because $C_1(\mathbf s, \mathbf X) =  \bar f_{b_2, \cdots, b_n}(\mathbf s)$ on $R_{b_2, \cdots, b_n}$). This is indeed true because if $b'_i \neq b_i$, the corresponding term would be negative in $f_{b'_2, \cdots, b'_n}(\mathbf s)$.
Finally, note that $c_1(\cdot) = \mathbb{E}[C_1(\cdot, \mathbf X)]$ is also convex since expectation preserves convexity.
\qed
}
\end{pf}
\blue{
Proposition \ref{prop: convexity of g1} shows that the $l_1$-type objective function is convex. The following example shows that this may not be true for the objective function $c_2$.
\begin{exmp}
\label{ex: non-convexity of g2}
Consider the special case of $n=3$ and let $X_1, X_2 > 0$ be positive scalars (which can be seen as degenerate distributions). We show that the function $c_2(s_2, s_3) := (X_1 - s_2)^2 + (\max \{X_1, s_2\} + X_2 - s_3)^2$ is not convex. Let $t = 0.5$ and $\mathbf s ^1 = (X_1 - \gamma, X_1 + 10X_2)$, $\mathbf s ^2 = (X_1, X_1 + 10X_2)$ and $\mathbf s ^3 = (X_1 + \gamma, X_1 + 10X_2)$ for some $0 < \gamma < \min \{X_1, 6X_2\}$. Substituting these values, we observe that $c_2(\mathbf s ^2) = c_2(t\mathbf s ^1 + (1 - t)\mathbf s ^3) > tc_2(\mathbf s ^1) + (1 - t)c_2(\mathbf s ^3)$.
\end{exmp}
}
\section{Useful Theorems and Propositions}
\label{sec:appendixb}
\begin{thm}[Theorem 2.3 of \cite{ArWe94}] Let $F: S \times \Xi \to (-\infty, \infty]$ be a measurable function and $P(d\xi)$ be a probability measure over the space $\Xi$ of random elements. We assume that $S$ is a metric space. Define $f(s) := \mathbb{E}[F(s, \xi)] = \int F(s, \xi)P(d\xi)$ and let $\xi_1, \cdots, \xi_m$ be independent samples of $\Xi$ drawn according to $P$. Suppose (1) $F(\cdot \,, \xi)$ is lower semicontinuous for fixed $\xi \in \Xi$ and (2) for each $s_0 \in S$ there exists an open set $N_0 \subseteq S$ and an integrable function $g_0: \Xi \to (-\infty, \infty)$ such that the inequality
	\begin{align*}
	F(s, \xi) \geq g(\xi)
	\end{align*}
	holds for all $s \in N_0$. Then, $\frac{1}{m}\sum_{j=1}^mF(\cdot, \xi_j)$ almost surely epi-converges to $f(\cdot)$.
\end{thm}
\begin{prop}[Proposition 7.26 of \cite{ShDeRu09}]
	Let $f_m, f: S \to (-\infty, \infty]$ where $S \subseteq \mathbb{R}^n$. Suppose that $f_m(\cdot)$ epi-converges to $f(\cdot)$. Then,
	\begin{align*}
	\limsup_m[\inf_s f_m(s)] \leq \inf_s f(s).
	\end{align*}
	Suppose further that (1) for some $\epsilon_m \downarrow 0$ there exists an $\epsilon_m-$minimizer $s_m$ of $f_m(\cdot)$ such that the sequence $s_m$ converges to a point $\bar{s}$. Then, $\bar{s} \in \text{argmin} f$ and
	\begin{align*}
	\lim_{m\to \infty}[\inf_s f_m(s)] = \inf_s f(s)
	\end{align*}
\end{prop}
\begin{prop}[Proposition 2.3.2 of \cite{BeNeOz03}]
	Let $S$ be a closed convex subset of $\mathbb{R}^n$, and let $f: \mathbb{R}^n \to (-\infty, \infty]$ be a closed convex function such that $f(s) < \infty$ for some $s \in S$. The set of minimizing points of $f$ over $S$ is nonempty and compact if and only if $S$ and $f$ have no common nonzero direction of recession.
\end{prop}

\end{document}